\documentclass[reqno,12pt,A4paper]{amsart}

\usepackage[english]{babel}
\usepackage{amsfonts, amsmath, amssymb, amsthm}
\usepackage{enumitem}

\usepackage{hyperref}

\usepackage{amsthm}
\newtheorem{definition}{Definition}%[section]
\newtheorem{proposition}[definition]{Proposition}%[section]
\newtheorem{theorem}[definition]{Theorem}%[section]
\theoremstyle{remark}

\numberwithin{equation}{section}

\def\R{\mathbb{R}}

\def\T{\mathrm{T}}
\def\End{\mathrm{End}}
\def\eCl{\mathrm{Cl}^0}

\def\GL{\mathrm{GL}}
\def\CO{\mathrm{CO}}
\def\PSO{\mathrm{PSO}}
\def\SO{\mathrm{SO}}
\def\Spin{\mathrm{Spin}}
\def\d{\mathrm{d}}
\def\L{\mathrm{L}}

\def\n{\nabla}

\title[Clifford-Weyl Structures]
{Local Geometry of Even Clifford Structures on Conformal Manifolds}

\author{Charles Hadfield}
\email{charles.hadfield@ens.fr}
\address{\'Ecole Normale Sup\'erieure, 45 rue d'Ulm,
75230 Paris cedex 05, France}
\author{Andrei Moroianu}
\email{andrei.moroianu@math.cnrs.fr}
\address{Laboratoire de Math\'ematiques de Versailles, UVSQ, CNRS, Universit\'e, 
Paris-Saclay, 78035 Versailles, France}

\begin{document}

\begin{abstract}
We introduce the concept of a Clifford-Weyl structure on a conformal manifold, which consists of an even Clifford structure parallel with respect to the tensor product of a metric connection on the Clifford bundle and a Weyl structure on the manifold. We show that the Weyl structure is necessarily closed except for some ``generic" low-dimensional instances, where explicit examples of non-closed Clifford-Weyl structures can be constructed. 

\bigskip


\smallskip

\end{abstract}

\maketitle

%%%%% %%%%% %%%%% %%%%% %%%%% %%%%% %%%%% %%%%% %%%%% %%%%% %%%%% %%%%%

\section{Introduction}

Even Clifford structures on Riemannian manifolds were introduced in \cite{MS} as a natural generalisation of almost Hermitian and almost quaternion Hermitian geometries and consist of a (locally defined) Euclidean vector bundle $(E,h)$ over $(M,g)$ together with an algebra bundle morphism $\varphi:\eCl(E,h)\to\End(\T M)$ mapping $\Lambda^2E$ into the bundle of skew-symmetric endomorphisms $\End^-(\T M)$. Alternatively, one can define even Clifford structures as structure group reductions, cf. \cite{ah} and \cite{agh} for details. Several articles appeared recently on this topic, see e.g. \cite{aha}, in which the twistor space is constructed, \cite{ahs} which studies even Clifford structures with large automorphism groups, \cite{gh} for rigidity and vanishing results, \cite{mp} for the classification of homogeneous even Clifford structures, or \cite{pp} and \cite{ppv} for relations with symmetric spaces and the Severi varieties.

The present paper studies even Clifford structures over conformal -- rather than Riemannian -- manifolds. We introduce the notion of a Clifford-Weyl structure, defined by the following data: a conformal manifold $(M,c)$, equipped with a Weyl structure $D$ and an even Clifford structure $(E,h,\varphi)$ which is parallel with respect to $D$ and some metric connection $\nabla^E$ on $E$. In the case where $D$ is exact (that is, the Levi-Civita connection of some metric in $c$), this notion is equivalent to that of a parallel even Clifford structure from \cite{MS}. Immediately the natural question to ask is under what conditions does this problem locally reduce to a problem in Riemannian geometry, i.e. under what conditions is $D$ closed? We show that there are six instances (called generic) where the presence of a Clifford-Weyl structure need not force the Weyl connection to be closed, and that in all other cases (called non-generic), the associated Weyl structure of a Clifford-Weyl structure is automatically closed.  More precisely, our results can be stated as follows:

\begin{theorem}\label{thm:theorem}
Suppose a conformal manifold of dimension $n$ carries a rank $r\ge 2$ Clifford-Weyl structure such that $(n,r)$ is different from $(2,2)$, $(4,2)$, $(4,3)$, $(4,4)$ and $(8,8)$. Then the associated Weyl connection is closed. The same conclusion holds if $(n,r)=(8,4)$, provided that the restriction of the Clifford morphism $\varphi$ to $\Lambda^2 E$ is not injective.
\end{theorem}

The cases excluded by this theorem are somehow generic, and they are treated in the following 

\begin{proposition}\label{prop:generic-cases}
\begin{enumerate}[leftmargin=*, label=(\roman*)]
\item Let $D$ be a Weyl structure on an oriented conformal manifold $(M,c)$ of dimension $2$, $4$ or $8$. Then $(M,c)$ carries a Clifford-Weyl structure of rank $r=2$ for $n=2$, $r=3$ or $r=4$ for $n=4$ and $r=8$ for $n=8$, whose associated Weyl structure is $D$.

\item Let $D$ be a Weyl structure on a conformal manifold $(M^n,c)$. Then there exists a Clifford-Weyl structure of rank $2$ on $(M,c)$ with associated Weyl structure $D$ if and only if $D$ preserves a complex structure compatible with $c$. If $n=4$, every complex structure $J$ compatible with $c$ is preserved by a unique Weyl structure $D^J$, which is closed if and only if $J$ is locally conformally K\"ahler.

\item Let $D$ be a Weyl structure on a conformal manifold $(M^8,c)$. Then there exists a Clifford-Weyl structure of rank $4$ whose Clifford morphism $\varphi:\eCl(E,h)\to\End(\T M)$ is injective upon restriction to $\Lambda^2 E$, if and only if $D$ is the adapted Weyl structure of a conformal product structure on $(M,c)$ with $4$-dimensional factors (cf. \cite{BMproducts}). 
\end{enumerate} 
\end{proposition}

The paper is organised as follows.
Section~\ref{sec:preliminaries} recalls several notions of even Clifford structures, defines Clifford-Weyl structures and introduces the required differential and algebraic objects from differential and conformal geometry. It finishes with a toy problem from Hermitian geometry which inspires the beginning of the proof of the theorem (and provides the proof for the case $n>4$, $r=2$). Section~\ref{sec:largerank} establishes the theorem for large rank $r\ge 5$ structures and Section~\ref{sec:smallrank} establishes the theorem in the remaining low rank setting. Section~\ref{sec:generic-cases} considers the generic cases of Proposition~\ref{prop:generic-cases} and shows that in each of these cases there are examples of Clifford-Weyl structures with non-closed associated Weyl structures.

%%%%% %%%%% %%%%% %%%%% %%%%% %%%%% %%%%% %%%%% %%%%% %%%%% %%%%% %%%%%
%%%%% %%%%% %%%%% %%%%% %%%%% %%%%% %%%%% %%%%% %%%%% %%%%% %%%%% %%%%%
%%%%% %%%%% %%%%% %%%%% %%%%% %%%%% %%%%% %%%%% %%%%% %%%%% %%%%% %%%%%

\section{Preliminaries}\label{sec:preliminaries}

Let $M$ be a smooth manifold. It is well-known that there is a one-to-one correspondence between isomorphism classes of oriented Euclidean rank $k$ vector bundles $(E,h)$ over $M$ and isomorphism classes of principal $\SO(k)$-bundles $P$ over $M$. By this correspondence, one may identify metric covariant derivatives on $(E,h)$ with connections on $P$. We denote 
\begin{align*}
\PSO(k)
:=
\begin{cases} 
			\SO(k)				& \hbox{if $k$ is odd}\\
			\SO(k)/\{\pm {\mathrm I}_k\}		&\hbox{if $k$ is even}\end{cases}
\end{align*}
and by a slight abuse of language we set the following

\begin{definition}
A \emph{locally defined} oriented Euclidean rank $k$ vector bundle over $M$ is a principal $\PSO(k)$-bundle over $M$.
\end{definition}

The terminology is justified by the fact that the structure group of a principal $\PSO(k)$-bundle can be reduced to $\SO(k)$ over any contractible open neighborhood $U$ of $M$, and thus gives rise to an oriented Euclidean rank $k$ vector bundle over $U$. 

If $E$ is a locally defined oriented Euclidean rank $k$ vector bundle over $M$ and $\rho:\PSO(k)\to \SO(N)$ is a group morphism, one obtains a rank $N$ oriented Euclidean vector bundle $\rho(E)$ over $M$ by enlarging the structure group of $E$ to $\SO(N)$ and considering the associated vector bundle. In particular, since the even-dimensional tensor powers of the standard representation of $\SO(k)$ on $\mathbb{R}^k$ descend to $\PSO(k)$, the even tensor powers of a locally defined oriented Euclidean vector bundle are globally defined vector bundles. 

We can now recall the definition of even Clifford structures on Riemannian manifolds, which have been introduced in \cite{MS}.

\begin{definition}
A rank $r\ge 2$ \emph{even Clifford structure} on a Riemannian manifold $(M^n,g)$ is an oriented, locally defined, rank $r$ Euclidean bundle $(E,h)$ over $M$ together with a non-vanishing algebra bundle morphism, called \emph{Clifford morphism}, $\varphi:\eCl(E,h)\to\End(\T M)$ which maps $\Lambda^2 E\subset \eCl(E,h)$ to the bundle of skew-symmetric endomorphisms $\End^-(\T M)$
\end{definition}

\begin{definition}
An even Clifford structure $(M,g,E,h,\varphi)$, is called \emph{parallel} if there exists a metric connection $\n^E$ on $(E,h)$ such that $\varphi$ is connection preserving with respect to $\n^E$ and the Levi-Civita connection $\n$ of $(M,g)$.
\end{definition}

As $\End^-(\T M)$ is invariant under a conformal change of the metric $g$, the notion of an even Clifford structure extends directly to the setting of a conformal manifold $(M,c)$. The condition of parallelism is transferred by considering Weyl connections giving what we term Clifford-Weyl structures.

\begin{definition}
A rank $r\ge2$ \emph{Clifford-Weyl structure} on a conformal manifold $(M^n,c)$,  is a tuple $(E,h,\varphi,\n^E,D)$ where
\begin{itemize}
\item $(E,h)$ is an oriented locally defined rank $r$ Euclidean bundle;
\item $\varphi:\eCl(E,h)\to\End(\T M)$ is an algebra bundle morphism sending $\Lambda^2 E$ to $\End^-(\T M)$;
\item $\n^E$ is a metric connection on $E$;
\item $D$ is a Weyl connection on $(M,c)$,
\end{itemize}
such that $\varphi$, seen as a section of $\eCl(E,h)^*\otimes \End(\T M)$, is parallel with respect to $\n^E\otimes D$.\end{definition}

Let $(E,h,\varphi,\n^E,D)$ be a rank $r$ Clifford-Weyl structure on a conformal manifold $(M^n,c)$ and let $L$ denote the weight bundle of $M$ (the real line bundle associated with the principal bundle of frames via the representation $|\det|^{1/n}$ of $\GL(n;\R)$, cf. \cite[Section 2]{BMproducts}). 

Consider a metric $g$ in the conformal class $c$.
Associated with $g$ we have the Levi-Civita connection $\n$ as well as the gauge $\ell$ (a section of $L$) and Lee form $\theta$ defined respectively by
\[
c = g\otimes \ell^2,
\qquad
D\ell = \theta \otimes \ell.
\]
This gives the useful formula $Dg = - 2\theta\otimes g$. Independent of the choice of metric in the conformal class, we have the Faraday form $F=\d\theta$. The Weyl structure $D$ is closed if and only if $D$ is locally the Levi-Civita connection of a metric in the conformal class. This is equivalent to $F=0$.

Let $\xi_1,\dots,\xi_r$ be a local oriented orthonormal frame for $(E,h)$. We introduce the collection of connection coefficients and the curvature two-forms
\[
\n^E \xi_j =: \sum_i\eta_{ij} \otimes \xi_i,
\quad
R^E \xi_j =: \sum_i \omega_{ij} \otimes \xi_i,
\] 
(with $\omega_{ij}=\d\eta_{ij}+\sum_k\eta_{ik}\wedge\eta_{kj}$),
and define endomorphisms $J_{ij} := \varphi(\xi_i\cdot \xi_j)$ where $\cdot$ denotes Clifford multiplication. Because of the relations in the Clifford algebra, and the fact that $\varphi$ is an algebra bundle morphism mapping $\Lambda^2E$ into $\End^-( \T M)$, the endomorphisms $J_{ij}$ are locally defined almost Hermitian structures on $M$ for $i\neq j$. Moreover, for mutually distinct indices $i,j,k$ we have 
\[(\xi_i\cdot \xi_j)\cdot (\xi_i\cdot \xi_k)=\xi_j\cdot\xi_k=-(\xi_i\cdot \xi_k)\cdot (\xi_i\cdot \xi_j)\]
and thus $J_ {ij}$ anticommutes with $J_{i_k}$. In particular, this shows that 
\begin{equation}\label{ac}
\langle J_{ij}, J_{ik}\rangle=0,\qquad\forall\ i\neq j\neq k\neq i,
\end{equation}
where $\langle \cdot, \cdot\rangle$ denotes, as usual, minus the trace of the product of two endomorphisms.

Let $e_1,\dots,e_n$ denote a local orthonormal frame for $(\T M,g)$. For each $J_{ij}$ with $i\neq j$, we obtain an associated non-degenerate two-form $\Omega_{ij}$:
\[
\Omega_{ij}(\cdot,\cdot):=g(J_{ij}\,\cdot,\cdot).
\]
(In the process of establishing \eqref{eq:beast}, as well in the the rank $r=4$ case, the calculations are simplified by summing indiscriminately over subscripts, for this we define $\Omega_{ii}:=0$.) Using the natural scalar product on the bundle of exterior forms induced by $g$, we obtain the Lefschetz-type operators for $i\neq j$,
\[
\L_{ij} := \Omega_{ij}\wedge,
\quad
\Lambda_{ij} := \L_{ij}^*=\tfrac12 \sum_{a=1}^n J_{ij}(e_a)\,\lrcorner\, e_a \,\lrcorner.
\]
For later use, notice that by the usual identification using the metric $g$ of $\Lambda^2(\T^*M)$ with $\End^-(\T M)$, the Lefschetz operator $\Lambda_{ij}$ acting on $\Lambda^2(\T^*M)$ (with $i\neq j$) is identified with $\tfrac12\langle J_{ij} , \cdot \rangle $ acting on $\End^-(\T M)$.

We finish this section with a toy problem: That when $r=2$ and $n>4$, the Weyl connection is closed. This is a standard fact in Hermitian geometry, but the proof below contains, at embryonic state, the main ideas of the proof of Theorem \ref{thm:theorem}. 

We choose as before a metric $g$ in the conformal class and drop the superfluous subscripts on $J_{12}$ and $\Omega_{12}$. First, the fact that $\varphi$ is $\n^E\otimes D$-parallel implies
\[
D J = D(\varphi(\xi_1\cdot \xi_2))= \varphi ( \n^E(\xi_1\cdot \xi_2) )
\]
and as
\[
\n^E (\xi_1\cdot \xi_2 )
= (\eta_{21}\otimes \xi_2)\cdot \xi_2 + \xi_1\cdot(\eta_{12}\otimes \xi_1)
= -(\eta_{12} + \eta_{21})1_{\eCl(E,h)} = 0,
\]
we conclude $J$ is parallel with respect to $D$.
Differentiating $\Omega$ with respect to the Weyl connection
\[
(D\Omega)(\cdot,\cdot) = (Dg)(J\cdot,\cdot) + g(DJ\cdot,\cdot)
\]
and using the formulae $DJ=0$ and $Dg = - 2\theta \otimes g$ gives $D\Omega=-2\theta\otimes \Omega$ which upon extracting the totally antisymmetric part yields
\[
\d \Omega = -2\theta\wedge\Omega.
\]
Differentiating this equation gives
\[
0=\d ^2\Omega = -2 F\wedge\Omega - 4\theta\wedge\theta\wedge\Omega
\]
and as $\Omega$ is non-degenerate with $n>4$, the equation $F\wedge\Omega=0$ forces $F=0$.

%%%%% %%%%% %%%%% %%%%% %%%%% %%%%% %%%%% %%%%% %%%%% %%%%% %%%%% %%%%%
%%%%% %%%%% %%%%% %%%%% %%%%% %%%%% %%%%% %%%%% %%%%% %%%%% %%%%% %%%%%
%%%%% %%%%% %%%%% %%%%% %%%%% %%%%% %%%%% %%%%% %%%%% %%%%% %%%%% %%%%%

\section{Large rank Clifford-Weyl structures}\label{sec:largerank}

For this section suppose that $(E,h,\varphi,\n^E,D)$ is a rank $r$ Clifford-Weyl structure on a conformal manifold $(M^n,c)$, with $r\ge 5$. The structure of $\eCl_r$ forces the dimension of the manifold to be a multiple of $8$. Apart from the generic situation $n=8$, $r=8$, which will be treated later on, we will show that the Weyl connection is closed. 

As in the previous section, we consider an arbitrary Riemannian metric $g$ in the conformal class $c$. Then $(E,h)$ is an even Clifford structure, and we may build Lefschetz-type operators as well as identify $\Lambda^2(\T^*M)$ with $\End^-(\T M)$.

The connection coefficients give for all $i\neq j$
\begin{align*}
\n^E (\xi_i\cdot \xi_j )
&= \sum_k \eta_{ki}\otimes (\xi_k\cdot \xi_j) + \eta_{kj}\otimes( \xi_i\cdot \xi_k) \\
&=\sum_{k\neq i,j} \eta_{ki} \otimes(\xi_k\cdot \xi_j )-\eta_{kj} \otimes( \xi_k\cdot \xi_i)
\end{align*}
and as the even Clifford structure is parallel,
\[
D J_{ij} = \sum_{k\neq i,j} \eta_{ki} \otimes J_{kj} - \eta_{kj}\otimes J_{ki}.
\]
Differentiating $\Omega_{ij}$ with respect to the Weyl connection
\[
(D\Omega_{ij})(\cdot,\cdot) = (Dg)(J_{ij}\cdot,\cdot) + g(DJ_{ij}\cdot,\cdot)
\]
and using the previous formula as well as the fundamental formula $Dg = - 2\theta \otimes g$ provides
\[
D\Omega_{ij} = - 2\theta\otimes \Omega_{ij} + \sum_{k\neq i,j} \eta_{ki}\otimes\Omega_{kj} - \eta_{kj}\otimes\Omega_{ki}.
\]
Taking the totally antisymmetric part of this equation (and recalling $\Omega_{ii}:=0$) gives
\[
\d\Omega_{ij} = - 2\theta\wedge \Omega_{ij} + \sum_k \eta_{ki}\wedge\Omega_{kj} - \eta_{kj}\wedge\Omega_{ki}.
\]
Differentiating this equation and replacing appearances of $\d\Omega_{ij}$ (as well as $\d\Omega_{kj}$ and $\d\Omega_{ki}$) using this same equation yields
\begin{align*}
2F\wedge \Omega_{ij} 
&=\sum_k \Big(
		2\theta\wedge \eta_{ki}\wedge \Omega_{kj} 
		+\d \eta_{ki} \wedge \Omega_{kj} - \eta_{ki} \wedge \d  \Omega_{kj}
		\Big)
	- \{i\leftrightarrow j\}
\\
&=\sum_k \Big(
		\d \eta_{ki} \wedge \Omega_{kj} - \eta_{ki} \wedge
			\sum_\ell \Big(
							\eta_{\ell k}\wedge \Omega_{\ell j} - \eta_{\ell j}\wedge\Omega_{\ell k}
			\Big)
		\Big)
	- \{i\leftrightarrow j\}
\\
&=\sum_k \Big(
		\d \eta_{ki} + \sum_\ell \eta_{k\ell}\wedge \eta_{\ell i}
		\Big)\wedge \Omega_{kj}
	- \{i\leftrightarrow j\}
\end{align*}
where $\{i\leftrightarrow j\}$ corresponds to the previously displayed term with indices $i$ and $j$ interchanged. The previous equation simplifies upon introducing the curvature two-forms $\omega_{ij}$ of $\n^E$ into
\[
2 F \wedge \Omega_{ij}
=\sum_k \omega_{ki}\wedge \Omega_{kj} - \omega_{kj}\wedge \Omega_{ki}.
\]
Writing this using the Lefschetz-type operators establishes, for all $i\neq j$,
\begin{align}\label{eq:LF}
2 \L_{ij} F 
=\sum_{k} \L_{kj}\omega_{ki} -\L_{ki}\omega_{kj}.
\end{align}

Assuming $i\neq j$, we apply $\Lambda_{ij}$ to \eqref{eq:LF}. Calculating $\Lambda_{ij}$ applied to the left hand side of \eqref{eq:LF} is aided by the $\mathfrak{sl}(2)$ structure of the Lefschetz-type operators. Specifically $2[\Lambda_{ij},\L_{ij}]F = (n-4) F$ and $\L_{ij}(2\Lambda_{ij} F)$ identifies, as an endomorphism, with $\langle J_{ij}, F\rangle J_{ij}$. In order to calculate $\Lambda_{ij}$ applied to the right hand side of \eqref{eq:LF}, we note that the sum in \eqref{eq:LF} may be taken over $k\neq i,j$ and we write 
\[
\Lambda_{ij}
=\tfrac12 \sum_a J_{ij}(e_a)\,\lrcorner\, e_a \,\lrcorner
\]
(recall that $e_1,\dots e_n$ denotes a local orthonormal frame for $(\T M,g)$). This gives, for $k\neq i,j$,
\begin{align*}
2\Lambda_{ij}\L_{kj}\omega_{ki}
&=  \sum_a J_{ij}(e_a) \,\lrcorner\, 
		\Big( \Omega_{kj}(e_a)\wedge \omega_{ki} + \Omega_{kj}\wedge \omega_{ki}(e_a)
		\Big)
\end{align*}
The summands in the previous display consist of four terms, the first of which vanishes because of \eqref{ac}. Developing the remaining three terms gives
\begin{align*}
2\Lambda_{ij}\L_{kj}\omega_{ki}
&=\sum_a  - \Omega_{kj}(e_a)\wedge (\omega_{ki}\circ J_{ij})(e_a) + (\Omega_{kj}\circ J_{ij})(e_a) \wedge \omega_{ki}(e_a) +  \langle J_{ij}, \omega_{ki} \rangle \Omega_{kj}.
\end{align*}
Testing against two tangent vectors shows
\begin{align*}
\sum_a  - \Omega_{kj}(e_a)\wedge (\omega_{ki}\circ J_{ij})(e_a) 
	&=\omega_{ki}(J_{ik}\cdot,\cdot) + \omega_{ki}(\cdot, J_{ik}\cdot)
	\\
\sum_a (\Omega_{kj}\circ J_{ij})(e_a) \wedge \omega_{ki}(e_a)
	&=\omega_{ki}(J_{ik}\cdot,\cdot) + \omega_{ki}(\cdot, J_{ik}\cdot)
\end{align*}
which, viewed as endomorphisms via the metric, are each precisely $[J_{ki},\omega_{ki}]$. Therefore, for $i,j,k$ distinct,
\[
\Lambda_{ij}\L_{kj}\omega_{ki}
=
[J_{ki},\omega_{ki}] +\tfrac 12 \langle _{ij}, \omega_{ki} \rangle \Omega_{kj}
\]
which establishes, for $i\neq j$,
\begin{align}\label{eq:beast}
(n-4) F + \langle J_{ij}, F\rangle \Omega_{ij}
=
\sum_{k\neq i,j}
[J_{ki},\omega_{ki}] + [J_{kj},\omega_{kj}]
+ \tfrac 12\langle J_{ij}, \omega_{ki}\rangle \Omega_{kj}
- \tfrac 12 \langle J_{ij}, \omega_{kj}\rangle \Omega_{ki}.
\end{align}

Working with \eqref{eq:beast}, we apply $2\Lambda_{ij}, 2\Lambda_{ia}, 2\Lambda_{ab}$ for $a,b$ different from $i,j$. For $2\Lambda_{ij}$ applied to \eqref{eq:beast} we remark that $\langle J_{ij}, J_{ij} \rangle=n$ while $J_{ij}$ is orthogonal to $J_{kj}$ and $J_{ki}$ for $i,j,k$ distinct hence
\begin{align*}
(2n-4) \langle J_{ij}, F\rangle
&=
\sum_{k\neq i,j}
\langle J_{ij} , [J_{ki},\omega_{ki}] \rangle
+
\langle J_{ij} , [J_{kj},\omega_{kj}] \rangle \\
&=
\sum_{k\neq i,j}
\langle [J_{ij} , J_{ki} ] ,\omega_{ki} \rangle
+
\langle [J_{ij} , J_{kj} ],\omega_{kj} \rangle
\end{align*}
establishing
\begin{align}\label{eq:JijF}
(2n-4) \langle J_{ij} , F \rangle
&=
2\sum_{k\neq i,j} \langle J_{kj},\omega_{ki} \rangle - \langle J_{ki} , \omega_{kj} \rangle,
\end{align}
For $2\Lambda_{ia}$ applied to \eqref{eq:beast}, we remark that $\langle J_{ia}, J_{ki} \rangle=-n\delta_{ak}$ for $k\neq i,j$ and importantly, as $r\ge 5$, the terms involving $\langle J_{ia}, J_{kj}\rangle$ vanish \cite[Equation 2]{MS}. Therefore
\begin{align*}
(n-4) \langle J_{ia}, F\rangle
&=
\tfrac n2 \langle J_{ij}, \omega_{aj}\rangle 
+
\sum_{k\neq i,j}		%\Big(
\langle [J_{ia}, J_{ki}], \omega_{ki}\rangle
+ 
\langle [ J_{ia}, J_{kj} ] , \omega_{kj}\rangle			%\Big)
\end{align*}
where
\[
\sum_{k\neq i,j}
\langle [J_{ia}, J_{ki}], \omega_{ki}\rangle
=
\sum_{k\neq i,j}
\langle J_{ka} - J_{ak}, \omega_{ki}\rangle
=
2\sum_{k\neq j} \langle J_{ka}, \omega_{ki} \rangle
\]
and since $ J_{ia}$ commutes with $J_{kj}$  for $k\neq i,j$ except when $k=a$,
\[
\sum_{k\neq i,j}
\langle [ J_{ia}, J_{kj} ] , \omega_{kj}\rangle
=
\langle [J_{ia}, J_{aj} ], \omega_{aj} \rangle
=
-2 \langle J_{ij}, \omega_{aj} \rangle
\]
establishing
\begin{align}\label{eq:JiaF}
(n-4) \langle J_{ia}, F \rangle
&=
(\tfrac n2 -2 )\langle J_{ij}, \omega_{aj}\rangle + 2\sum_{k\neq j} \langle J_{ka}, \omega_{ki} \rangle,
\end{align}
For $2\Lambda_{ab}$ applied to \eqref{eq:beast}, we again use the large rank hypothesis $r\ge 5$. Indeed, due to this condition, for $k\neq i,j$, terms involving $\langle J_{ab},J_{kj}\rangle$ and $\langle J_{ab}, J_{ki}\rangle$ vanish. So
\begin{align*}
(n-4) \langle J_{ab}, F \rangle
&=
\sum_{k\neq i,j}
\langle [J_{ab},J_{ki}],\omega_{ki}\rangle + \langle [J_{ab},J_{kj}],\omega_{kj}\rangle \\
&= 
\sum_{k\in\{a,b\}}
\langle [J_{ab},J_{ki}],\omega_{ki}\rangle + \langle [J_{ab},J_{kj}],\omega_{kj}\rangle
\end{align*}
and developing the four terms from the summation in the preceding display establishes
\begin{align}\label{eq:JabF}
(n-4) \langle J_{ab}, F \rangle
&=
2\langle J_{bi},\omega_{ai}\rangle - 2\langle J_{ai},\omega_{bi}\rangle + 2\langle J_{bj},\omega_{aj}\rangle - 2\langle J_{aj},\omega_{bj}\rangle.
\end{align}

Armed with the preceding numbered equations, we may establish the orthogonality between $J_{ij}$ and $F$. From \eqref{eq:JabF}, by collecting the first two terms, and collecting the second two terms, we see that $\langle J_{bi},\omega_{ai}\rangle - \langle J_{ai},\omega_{bi}\rangle$ is independent of $i$ so
\[
(n-4) \langle J_{ab}, F \rangle
=
4\langle J_{bi},\omega_{ai}\rangle - 4\langle J_{ai},\omega_{bi}\rangle
\]
Summing the previous display over $i \neq a,b$ and changing the notation of indices $a,b,i\to i,j,k$ gives
\[
(r-2)(n-4) \langle J_{ij}, F \rangle 
=
4 \sum_{k\neq i,j} \langle J_{kj}, \omega_{ki} \rangle - \langle J_{ki},\omega_{kj} \rangle.
\]
Comparing this equation with \eqref{eq:JijF} provides the constraint
\[
(r-2)(n-4) \langle J_{ij}, F \rangle = 2(2n-4) \langle J_{ij}, F\rangle
\]
whence $\langle J_{ij}, F\rangle=0$ unless $4(2n-4)=2(r-2)(n-4)$. As $r\ge 5$ and $n$ is a multiple of $8$, the only obstructive case is the generic case $n=8$, $r=8$, which was excluded. Therefore 
\begin{align}
\langle J_{ij}, F\rangle=0\qquad\forall\ i\neq j.
\end{align}
Updating \eqref{eq:JiaF} and \eqref{eq:JabF} using this orthogonality, we obtain a pair symmetry from $\eqref{eq:JabF}$
\begin{align}
\label{eq:pairsymmetry}
\langle J_{ia},\omega_{ja}\rangle = \langle J_{ja},\omega_{ia}\rangle,
\qquad
\forall\ \textrm{$i,j,a$ distinct.}
\end{align}
which, upon switching the variables $j,a$ in \eqref{eq:JiaF}, provides
\[
(2-\tfrac n2) \langle J_{ia}, \omega_{ja} \rangle = -2 \langle J_{aj}, \omega_{ai} \rangle + 2 \sum_k \langle J_{kj}, \omega_{ki} \rangle,
\]
giving
\begin{align}
\label{eq:JiaFupdated}
(2-\tfrac n4) \langle J_{ia}, \omega_{ja} \rangle = \sum_k \langle J_{ik}, \omega_{jk} \rangle,
\qquad\forall\ 
\textrm{$i,j,a$ distinct.}
\end{align}
Therefore if $n=8$, the sum on the right hand side vanishes, while if $n\neq 8$, $\langle J_{ia}, \omega_{ja} \rangle$ is independent of $a\neq i,j$ hence
\[
(2-\tfrac n4) \langle J_{ia}, \omega_{ja} \rangle = (r-2) \langle J_{ia}, \omega_{ja} \rangle
\]
and so $J_{ia}$ is orthogonal to $\omega_{ja}$. It thus turns out that the sum $\sum_k \langle J_{ik}, \omega_{jk} \rangle$ vanishes no matter what the dimension $n$ is.

As a penultimate result, we remark that upon summation over $j\neq i$ (for $i$ fixed), the final two terms of \eqref{eq:beast} vanish:
\begin{align}\label{eq:last2vanish}
\sum_{j\neq i}\left(
			\sum_{k\neq i,j}
				\langle J_{ij}, \omega_{ki}\rangle J_{kj}
				-\langle J_{ij}, \omega_{kj}\rangle J_{ki}
		\right)
=0.
\end{align}
In fact, the previous display naturally splits into two collections of summations, each collection vanishing independently as we now show. 
The first collection of summations in \eqref{eq:last2vanish} may be written as the sum over $j,k$ both different from $i$ and from each-other:
\[
\sum_{j\neq i}\left(\sum_{k\neq i,j}
\langle J_{ij}, \omega_{ki}\rangle J_{kj}\right)
=
\sum_{\substack{j,k\neq i \\ j\neq k}}
\langle J_{ij}, \omega_{ki}\rangle J_{kj}
\]
which thus vanishes as $\langle J_{ij}, \omega_{ki}\rangle$ is symmetric in $j,k$ due to \eqref{eq:pairsymmetry} while $J_{jk}$ is antisymmetric in $j,k$. Considering the second collection of summations in \eqref{eq:last2vanish}, we rearrange the summation,
\begin{align*}
\sum_{j\neq i} \left( \sum_{k\neq i,j}
				\langle J_{ij}, \omega_{kj}\rangle J_{ki}
			\right)
&=
\sum_{j\neq i} \left( \sum_{k\neq i}
				\langle J_{ij}, \omega_{kj}\rangle J_{ki}
			\right)\\
&=
\sum_{k\neq i} \left( \sum_{j\neq i}
				\langle J_{ij}, \omega_{kj}\rangle 
			\right) J_{ki}\\
&=
\sum_{k\neq i,j} \left( \sum_j
				\langle J_{ij}, \omega_{kj}\rangle 
			\right) J_{ki}
\end{align*}
and by the remark following \eqref{eq:JiaFupdated}, the preceding display vanishes and provides \eqref{eq:last2vanish}.

We may now establish the result. By defining $S_i:=\sum_k [J_{ki},\omega_{ki}]$, \eqref{eq:beast} now reads
\[
(n-4)F = S_i + S_j - 2[J_{ij}, \omega_{ij} ] + \tfrac 12 \sum_{k\neq i,j} \langle J_{ij}, \omega_{ki} \rangle J_{kj} - \langle J_{ij}, \omega_{kj} \rangle J_{ki}.
\]
Keeping $i$ fixed and summing over $j\neq i$, making use of \eqref{eq:last2vanish}, we obtain
\begin{align*}
(r-1)(n-4) F
&= \sum_{j\neq i} \left( S_i + S_j - 2[J_{ij}, \omega_{ij}] \right)
\\
&= (r-4)S_i + \sum_{j} S_j
\end{align*}
which implies (as $r\neq 4$) that $S_i$ is independent of $i$ and proportional to $F$:
\[
(r-1)(n-4) F = 2(r-2)S_i.
\]
Equation \eqref{eq:beast} thus develops to
\[
\frac{2(n-4)}{r-2} F = 4 [ J_{ij}, \omega_{ij} ] -  \sum_{k\neq i,j}
\langle J_{ij}, \omega_{ki}\rangle J_{kj}
-\langle J_{ij}, \omega_{kj}\rangle J_{ki}.
\]
Commuting $F$ with $J_{ij}$ we see that
\begin{align*}
\frac{2(n-4)}{r-2}  FJ_{ij} &= 4( J_{ij} \omega_{ij} J_{ij} + \omega_{ij}) - \sum_{k\neq i,j}
\langle J_{ij}, \omega_{ki}\rangle J_{ki}
+\langle J_{ij}, \omega_{kj}\rangle J_{kj},
\\
\frac{2(n-4)}{r-2}  J_{ij} F &=-4(\omega_{ij} + J_{ij}\omega_{ij} J_{ij}) + \sum_{k\neq i,j}
\langle J_{ij}, \omega_{ki}\rangle J_{ki}
+\langle J_{ij}, \omega_{kj}\rangle J_{kj}.
\end{align*}
Therefore $F$ anticommutes with $J_{ij}$ for every $i\neq j$. By taking some $k$ different from both $i$ and $j$ we get that $F$ commutes with $J_{ik}J_{jk}=J_{ij}$. Hence $F=0$, thus proving Theorem \ref{thm:theorem} when the rank of the Clifford-Weyl structure is at least $5$.

\section{Low rank Clifford-Weyl structures}\label{sec:smallrank}

We consider now the remaining cases from Theorem \ref{thm:theorem}. If the rank of the Clifford-Weyl structure is $2$, this follows directly from our toy problem described at the end of Section \ref{sec:preliminaries}.

That $D$ is closed for $r=3$ and $n\ge 8$ is a standard result in quaternion Hermitian Weyl geometry (or locally conformally quaternion K\"ahler geometry) \cite{Ornea}. We present a proof which can also be adapted to the case $r=4$. 

Define $\Omega:=\Omega_{12}^2 + \Omega_{23}^2 + \Omega_{31}^2$ to be the fundamental four-form (or Kraines form) of quaternion Hermitian geometry. By \eqref{eq:LF}, which continues to hold for $r=3$, we obtain
\begin{align*}
\Omega_{12}^2\wedge F = \tfrac 12\Omega_{12}\wedge\Omega_{13}\wedge\omega_{32} - \tfrac 12\Omega_{12}\wedge\Omega_{23}\wedge\omega_{31}
\end{align*}
and cyclically commuting $(1,2,3)$ gives two similar equations. Upon summation, cancellations give $\Omega\wedge F=0$ and, as the fundamental four-form is well-known to be non-degenerate (and $n\neq 4$), $F=0$. (Alternatively, if one follows the derivation of \eqref{eq:LF}, one obtains similar equations for $\d \Omega_{ij}$ in terms of the connection coefficients $\eta_{ij}$ which result in the equation $\d \Omega = -4\theta\wedge \Omega$. Differentiating a second time gives $\Omega\wedge F=0$.)

Finally, if $(E,h,\nabla^E,\varphi)$ is a rank 4 Clifford-Weyl structure and $n\ge 8$, let us consider $A\in\End(\T M)$ to be the image under $\varphi$ of the volume element of $E$. From the properties of $\varphi$, $A$ is a symmetric involution, hence the tangent bundle splits into a direct sum $\T M=T^+\oplus T^-$ of the $\pm 1$ eigenspaces of $A$. If either $T^+$ or $T^-$ are of dimension zero, then the rank 4 even Clifford structure is effectively a rank 3 even Clifford structure and the result follows from the previous paragraph. We may thus assume the decomposition of $\T M$ is non-trivial. In particular $\varphi$ is injective upon restriction to $\Lambda^2 E$, so we only need to consider the case $n\ge 12$. 

We construct quaternionic structures on $T^\pm$
\[
J_{12}^\pm = \mp \tfrac 12 ( J_{14} \pm J_{23} ),
\quad
J_{31}^\pm = \mp \tfrac 12 ( J_{13} \mp J_{24} ),
\quad
J_{23}^\pm = \mp \tfrac 12 ( J_{12} \pm J_{34} )
\]
which vanish upon restriction to $T^\mp $. We may thus define two four-forms $\Omega_\pm\in \Lambda^4(T^\pm)^*$ as in the case of quaternion Hermitian geometry and set $\Omega=\Omega_+ + \Omega_-$. We decompose the exterior algebra $\Lambda^* M = \oplus (\Lambda^p (T^+)^* \oplus \Lambda^q (T^- )^*)$ and say that elements of $\Lambda^p (T^+)^* \oplus \Lambda^q (T^-)^*$ are of type $(p,q)$. The decomposition of $\Lambda^2M$ enables us to write $F=F_+ + F_m + F_-$ where $F_+$, $F_m$, $F_-$ are respectively of type $(2,0)$, $(1,1)$, $(0,2)$. Using this we calculate $\Omega\wedge F$ (whose 6 pieces are of distinct type). Meanwhile, we remark that
\[
\Omega = \tfrac 12 \sum_{i<j} \Omega_{ij}^2 = \tfrac 14 \sum_{i,j} \Omega_{ij}^2
\]
and via \eqref{eq:LF}, which continues to hold for $r=4$,
\[
\Omega_{ij}^2\wedge F = \tfrac 12 \sum_{k} \Omega_{ij}\wedge\Omega_{ik}\wedge\omega_{kj} - \Omega_{ij}\wedge\Omega_{jk}\wedge\omega_{ki}
\]
Summing over $i,j$ the second term in the previous sum, $\sum_{i,j,k} \Omega_{ij}\wedge\Omega_{jk}\wedge\omega_{ki}$ may be written, under a permutation $i,j,k\to k,i,j$, as the first term $\sum_{i,j,k} \Omega_{ij}\wedge\Omega_{ik}\wedge\omega_{kj}$ hence
\[
\Omega\wedge F = 0.
\]
Since they have different types, each of the six terms in the expansion of $\Omega\wedge F$ also individually vanish. As $M$ is at least $12$-dimensional with both subbundles $T^\pm $ being non-trivial, we deduce from $\Omega_-\wedge F_+=0$ and $\Omega_+\wedge F_-=0$ that $F_\pm=0$. And as one of the subbundles $T^\pm $ has rank larger than $4$ (say $T^+$) then $F_m=0$ (from $\Omega_+\wedge F_m=0$). 

This finishes the proof of Theorem \ref{thm:theorem} when the rank of the Clifford-Weyl structure is $2$, $3$ or $4$.

%%%%% %%%%% %%%%% %%%%% %%%%% %%%%% %%%%% %%%%% %%%%% %%%%% %%%%% %%%%%
%%%%% %%%%% %%%%% %%%%% %%%%% %%%%% %%%%% %%%%% %%%%% %%%%% %%%%% %%%%%
%%%%% %%%%% %%%%% %%%%% %%%%% %%%%% %%%%% %%%%% %%%%% %%%%% %%%%% %%%%%

\section{Generic cases}\label{sec:generic-cases}

In this final section we prove Proposition~\ref{prop:generic-cases} and, in the process, show examples of Clifford-Weyl structures with non-closed associated Weyl covariant derivatives.

$(i)$ If $M$ has dimension $2$, we define $(E,h)$ to be the trivial rank $2$ Euclidean vector bundle with trivial flat connection $\nabla^E$, and $\varphi:\Lambda^2E\to \End^-(\T M)$ by $\varphi(\xi_1\wedge \xi_2):=J$, where $\xi_1,\xi_2$ is an oriented orthonormal frame of $E$ and $J$ is the rotation in $\T M$ by $\pi/2$ in the positive direction determined by $c$. Since $DJ=0$, $(E,h,\varphi,\nabla^E,D)$ is a rank $2$ Clifford-Weyl structure.

If $M$ has dimension $4$, we define $E:=\Lambda^+M\otimes L^2$ (the bundle of self-dual two-forms of conformal weight $0$), $\nabla^E$ to be the  covariant derivative induced by $D$ on $E$ and $h$ to be the canonical scalar product induced by $c$ on $E$. Since $\Lambda^2M\otimes L^2$ is canonically isomorphic to $\End^-(\T M)$, $E$ is in fact a rank $3$ sub-bundle of the bundle of skew-symmetric endomorphisms of $M$. Moreover, since $E$ is oriented, the metric $h$ provides an identification of $\Lambda^2E$ with $E$, and thus a map $\varphi:\Lambda^2E\to \End^-(\T M)$. It is straightforward to check that this map extends to an algebra morphism from $\eCl(E,h)$ to $\End(\T M)$ which is tautologically parallel with respect to $\nabla^E\otimes D$, thus defining a rank $3$ Clifford-Weyl structure. 

Moreover, every rank $3$ Clifford-Weyl structure $(E,h,\nabla^E,\varphi,D)$ determines in a tautological way a rank $4$ Clifford-Weyl structure $(\tilde E,\tilde h,\nabla^{\tilde E},\tilde\varphi,D)$ where $\tilde E=E\oplus \mathbb{R}$ with induced metric $\tilde h$ and connection $\nabla^{\tilde E}$, and $\tilde\varphi$ is defined on $\Lambda^2\tilde E\simeq \Lambda^2 E\oplus E$ by $\tilde\varphi=\varphi$ on $\Lambda^2 E$ and $\tilde\varphi=\varphi\circ\ast$ on $E$, where $\ast$ denotes the Hodge isomorphism $\ast:E\to\Lambda^2E$.

If $M$ has dimension $8$, we define $E=\Sigma^+_0M$ (the bundle of real half-spinors of conformal weight $0$, cf. \cite{BHMMM}) and $\nabla^E$ and $h$ to be the covariant derivative and the scalar product induced on $E$ by $D$ and $c$. Of course, if $M$ is not spin, $E$ is only locally defined, but $\Lambda^2E$ is always globally defined. We consider the map $\varphi:\Lambda^2E\to \End^-(\T M)$ defined by 
$$\varphi(\psi\wedge\phi):=X\mapsto -\sum_{i=1}^8 h(\ell^{-2}e_i\cdot X\cdot\psi,\phi)\, e_i-h(\psi,\phi)\, X,$$
where $\ell$ is a local section of $L$ and $e_i$ is a local frame of $\T M$ satisfying $c(e_i,e_j)=\ell^2 \delta_{ij}$. The map $\varphi$ is tautologically parallel with respect to $\nabla^E\otimes D$. Moreover, $\varphi$ extends to an algebra morphism from $\eCl(E,h)$ to $\End(\T M)$. Indeed, in order to check the universality property for the even Clifford algebra (\cite[Lemma A.1]{MS}), consider local sections $\psi$, $\phi$ and $\xi$ of $E$ such that $\psi$ is orthogonal to $\phi$ and $\xi$ and $h(\psi,\psi)=1$. Then $\{\ell^{-1}e_i\cdot \psi\}$ is a local orthonormal basis of the zero-weight half-spin bundle $\Sigma^-_0M$ (whose metric is also denoted by $h$) and thus
\begin{align*}
[\varphi(\psi\wedge\phi)\circ \varphi(\psi\wedge\xi)](X)
&=-\sum_{i=1}^8 h(\ell^{-2}e_i\cdot X\cdot\psi,\xi)\, \varphi(\psi\wedge\phi)(e_i)\\
&=\sum_{i,j=1}^8 h(\ell^{-2}e_i\cdot X\cdot\psi,\xi)h(\ell^{-2}e_j\cdot e_i\cdot\psi,\phi)\, e_j\\
&=\sum_{i,j=1}^8 h(\ell^{-2}X\cdot e_i\cdot\psi+2\ell^{-2}c(e_i,X)\psi,\xi)h(\ell^{-2}e_i\cdot\psi,e_j\cdot\phi)\, e_j\\
&=-\sum_{j=1}^8 h(\ell^{-2}X\cdot\xi,e_j\cdot\phi)\, e_j\\
&=-\varphi(\xi\wedge\phi)(X)-h(\phi,\xi)\, X\\
&=\varphi(\phi\wedge\xi)(X)-h(\phi,\xi)\, X.
\end{align*}
This shows that $(E,h,\varphi,\nabla^E,D)$ is a rank $8$ Clifford-Weyl structure on $M$.

$(ii)$ With any Clifford-Weyl structure of rank $2$ on $M$ one can associate the image, $J$, of the volume form of $E$ through the Clifford morphism $\varphi$. Clearly $J$ is an almost complex structure on $M$ compatible with $c$ and $D$-parallel. On the other hand, every almost complex structure preserved by a torsion-free connection is integrable.

Conversely, if $D$ is a Weyl structure on $(M,c)$ and $J$ is a $D$-parallel Hermitian structure, we define as before a rank $2$ Clifford-Weyl structure on $M$ by taking $(E,h)$ to be the trivial rank $2$ Euclidean vector bundle with trivial flat connection $\nabla^E$, and $\varphi:\Lambda^2E\to \End^-(\T M)$ defined by the fact that it maps the volume form of $E$ onto $J$.

For the second point, recall that on 4-dimensional conformal manifolds, every complex structure $J$ compatible with the conformal structure is preserved by a unique Weyl covariant derivative $D$ (see e.g. the proof of \cite[Lemma 5.7]{BMproducts}) which is closed if and only if $(J,c)$ is locally conformally K\"ahler.

$(iii)$ If $(E,h,\varphi,\nabla^E,D)$ is a Clifford-Weyl structure of rank $4$ with $\varphi$ injective on $\Lambda^2E$ on an $8$-dimensional conformal manifold $(M,c)$, then the image of the volume form of $E$ through $\varphi$ is a $D$-parallel involution of $\T M$ whose eigenbundles are $4$-dimensional $D$-parallel distributions. By \cite[Theorem 4.3]{BMproducts}, $(M,c)$ has a conformal product structure defined by these two distributions. 

Conversely, every conformal product structure on $(M,c)$ with 4-dimensional distributions $T^\pm $ defines a unique Weyl connexion $D$ (called the adapted Weyl structure in \cite[Definition 4.4]{BMproducts}) such that the splitting $\T M=T^+\oplus T^-$ is $D$-parallel. We obtain in this way a structure group reduction from $\CO(8)$ to $G:=\CO(8)\cap(\CO(4)\times \CO(4))$, that is, a $G$-principal bundle $P$ over $M$ and a connection induced by $D$ on $P$. Since $\CO(4)=\mathbb{R}^*\times (\Spin(3)\times_{\mathbb{Z}/2\mathbb{Z}} \Spin(3))$, the projections from $\mathbb{R}^*\times (\Spin(3)\times \Spin(3))$ to the second and third factors respectively define group morphisms $i_l$ and $i_r$ from $\CO(4)$ to $\SO(3)$. 
Let $E$ denote the (locally defined) rank $4$ Euclidean vector bundle over $M$ associated with $P$ via the group morphism $i_l\times i_r$ from $G$ to $\SO(3)\times \SO(3)=\PSO(4)$, and let $\nabla^E$ denote the induced covariant derivative. By construction, $\Lambda^2E$ is globally defined, and isomorphic to the weightless bundle $(\Lambda^+(T^+)\otimes L^{-2})\oplus(\Lambda^-(T^-)\otimes L^{-2})$. The composition of this isomorphism with the canonical inclusion in $\Lambda^+(\T M)\otimes L^{-2}=\End^-(\T M)$ yields as before a Clifford morphism, which is $\nabla^E\otimes D$-parallel by naturality of the construction.

Examples of conformal products with non-closed adapted Weyl structures can be easily constructed. Take $(M_1,g_1)$ and $(M_2,g_2)$ two $4$-dimensional Riemannian manifolds and let $M=M_1\times M_2$ with conformal class $c=[e^fg_1+g_2]$ where $f$ is any smooth map on $M$. Then the adapted Weyl structure of this conformal product structure is closed if and only if there exist functions $f_i$ on $M_i$ such that $f=\pi_1^*(f_1)+\pi_2^*(f_2)$ where $\pi_i:M\to M_i$ are the canonical projections (see \cite[Section 6.1]{BMproducts} for details).

%%%%% %%%%% %%%%% %%%%% %%%%% %%%%% %%%%% %%%%% %%%%% %%%%% %%%%% %%%%%
%%%%% %%%%% %%%%% %%%%% %%%%% %%%%% %%%%% %%%%% %%%%% %%%%% %%%%% %%%%%
%%%%% %%%%% %%%%% %%%%% %%%%% %%%%% %%%%% %%%%% %%%%% %%%%% %%%%% %%%%%

%%%%% %%%%% %%%%% %%%%% %%%%% %%%%% %%%%% %%%%% %%%%% %%%%% %%%%% %%%%%
%%%%% %%%%% %%%%% %%%%% %%%%% %%%%% %%%%% %%%%% %%%%% %%%%% %%%%% %%%%%
%%%%% %%%%% %%%%% %%%%% %%%%% %%%%% %%%%% %%%%% %%%%% %%%%% %%%%% %%%%%

%%%%% %%%%% %%%%% %%%%% %%%%% %%%%% %%%%% %%%%% %%%%% %%%%% %%%%% %%%%%

\def\arXiv#1{\href{http://arxiv.org/abs/#1}{arXiv:#1}}

\bibliographystyle{siam}

%%%%% %%%%% %%%%% %%%%% %%%%% %%%%% %%%%% %%%%% %%%%% %%%%% %%%%% %%%%%
%%%%% %%%%% %%%%% %%%%% %%%%% %%%%% %%%%% %%%%% %%%%% %%%%% %%%%% %%%%%
%%%%% %%%%% %%%%% %%%%% %%%%% %%%%% %%%%% %%%%% %%%%% %%%%% %%%%% %%%%%

\end{document}